\documentclass[12pt]{amsart}

\usepackage{amsmath,amssymb}
\usepackage[latin1]{inputenc}
\usepackage{url}
\newcommand{\urlBiBTeX}[1]{\url{#1}}

\newtheorem{thm}{Theorem}
\newtheorem{lem}[thm]{Lemma}
\newtheorem{prop}[thm]{Proposition}

\newcommand{\Z}{{\mathbf Z}}
\newcommand{\Q}{{\mathbf Q}}

\newcommand{\R}{{\mathbf R}}

\newcommand{\vol}{\operatorname{vol}}

\newcommand{\Sp}{ \Omega_p}
\newcommand{\Sa}{ \Omega_{q_1}}
\newcommand{\Sb}{ \Omega_{q_2}}

\newcommand{\Sq}{ \Omega_q}

\newcommand{\LL}{{\mathcal L}}

\newcommand{\AG}{{\mathcal AG}}

\newcommand{\vh}{{\bf h}}

\newcommand{\lattice}{L}
\newcommand{\lstar}{L^\times}

\newcommand{\eps}{\varepsilon}

\newcommand{\fp}{\mathbb{F}_p}
\newcommand{\kk}{\fp}
\newcommand{\kbar}{\overline{\mathbb{F}_p}}
\newcommand{\gal}{\operatorname{Gal}}
\newcommand{\fix}{\operatorname{Fix}}
\newcommand{\frob}{\operatorname{Frob}}
\newcommand{\m}{\mathfrak{m}}
\newcommand{\M}{\mathfrak{M}}

\newcommand{\cnot}{C_0}
\newcommand{\sym}{\tilde{R}}

\newcommand{\codim}{\operatorname{codim}}

\newcommand{\norm}{\operatorname{N}}
\newcommand{\cboundexponent}{{d \binom{k}{2} |\sym|}}

\newcommand{\rmkdd}[1]{}
\newcommand{\rmkd}[1]{}

\title{Poisson spacing statistics for value sets of polynomials}

\author{P\"ar Kurlberg}

\email{kurlberg@math.kth.se}
\address{Department of Mathematics\\ 
Royal Institute of Technology\\
SE-100 44 Stockholm  \\
Sweden}

\thanks{Author supported in part by 
the G\"oran Gustafsson Foundation, 
the National Science Foundation (DMS 0071503),
the Royal Swedish Academy of Sciences, 
  and the Swedish   Research Council.}

\date{November 25, 2005}

\begin{document}


\begin{abstract}
  If $f$ is a polynomial with integer coefficients and $q$ is an
  integer, we may regard $f$ as a map from $\Z/q\Z$ to $\Z/q\Z$.
  We show that the distribution of the (normalized) spacings between
  consecutive elements in the image of these maps becomes {\em
    Poissonian} as $q$ tends to infinity along any sequence of square
  free integers such that the 
mean spacing modulo $q$ tends to infinity.

\end{abstract}

\maketitle





\section{Introduction }

Let $f$ be a polynomial with integer coefficients.  Given an integer
$q$, we may regard $f$ as a map from $\Z/q\Z$ to $\Z/q\Z$, and the
image of this map will be denoted the {\em image of $f$ modulo $q$}.
The purpose of this paper is to investigate the distribution of
spacings between consecutive elements in the image of $f$ modulo $q$
as $q$ tends to infinity along square free integers.
The main emphasis will be placed on
the highly composite case, i.e., by letting $q$ tend to infinity in
such a way that the number of prime factors of $q$ also tends to
infinity. 

The case $f(x) = x^2$ 
and $q$ prime was investigated by Davenport. 
In
\cite{Davenport-quadratic,davenport-character-sums}
he proved 
%
%
that the probability of two consecutive squares being spaced
$h$ units apart tends to $2^{-h}$ as $q \to \infty$.  
We may interpret this
as if spacings between squares modulo prime $q$ behave like gaps between
heads in a sequence of fair coin flips.

The case $f(x) = x^2$ and $q$
 highly composite was studied by
Rudnick and the author  in \cite{squares1, squares2}.  If we let
$\omega(q)$ be the number of distinct prime factors of $q$, then the
number of squares modulo $q$ equals $\prod_{p|q}\frac{p+1}{2}$, and
the average spacing between the squares is given by 
$$
s_q =
\frac{q}{
\prod_{p|q}\frac{p+1}{2}}
=
2^{\omega(q)} \prod_{p|q}\frac{p}{p+1}.
$$
Hence $s_q \to \infty$  as $\omega(q) \to \infty$, so
we would expect that the probability of two squares being $1$ unit
apart vanishes as $\omega(q) \to \infty$, and it is thus natural to
normalize so that the mean spacing is one.  
A natural statistical model for the spacings is then given by looking
at random points in $\R/\Z$; for independent uniformly distributed
numbers in $\R/\Z$, the {\em normalized} spacings are said to be
Poissonian.  In particular, the distribution $P(s)$ of spacings
between consecutive points is that of a Poisson arrival process, i.e.,
$P(s) = e^{-s}$, and the joint distribution of $l$ consecutive
spacings is a product $l$ independent exponential random variables
(see \cite{feller-probability}).  Using Davenport's result together
with the heuristic that ``primes are independent'', it is seems
reasonable to expect that the distribution of the normalized spacings
between squares modulo $q$ becomes Poissonian in the limit $s_q \to
\infty$, and the main result of \cite{squares1} is that this is indeed
the case for squarefree $q$ (the general case is treated in
\cite{squares2}.) 

%
What can be said about more general polynomials 
$f \in \Z[x]$?  
For $p$  prime, let
$$
\Sp := \{ t \in  \Z/p\Z : t = f(x) \text{ for some $y \in  \Z/p\Z$} \}
$$
be the  image of $f$ modulo $p$.
Given $k \geq 2$ and integers $h_1, h_2, \ldots, h_{k-1}$, let
$$
N_k((h_1,h_2, \ldots, h_{k-1}), p) := 
|\{ t \in \Sp : t+h_1, \ldots, t+h_{k-1}  \in \Sp
\}|
$$
be the counting function for the number of $k$-tuples of elements in
the image of the form $t, t+h_1, \ldots, t+h_{k-1}$.
Letting $s_p := p/|\Sp|$ denote the average gap modulo $p$, the
``probability'' of an element being in the image is $1/s_p$.  Thus, if
the conditions $t \in \Sp, t+h_1 \in \Sp, \ldots, t+h_{k-1} \in \Sp$
are independent, we would expect $N_k((h_1,h_2, \ldots, h_{k-1}),
p)$ to be of size $p/s_p^k$, and a natural analogue of Davenport's 
result is then that 
\begin{equation}
\label{eq:Nk-bound}
N_k((h_1,h_2, \ldots, h_{k-1}),p) = p/s_p^k + o(p)  
\end{equation}
as $p\to \infty$  provided that $0, h_1, \ldots, h_{k-1}$ are
distinct modulo $p$.
%
In \cite{crtpoisson} Granville and the author proved that 
\begin{equation}
\label{eq:morse-mod-p-asymptotic}
N_k((h_1, h_2, \ldots, h_{k-1}),p) =
p/s_p^k + O_{f,k}(\sqrt{p})  
\end{equation}
holds if $f$ is a Morse polynomial and $0, h_1, \ldots, h_{k-1}$ are
distinct modulo $p$.
Using this, Poisson spacings for the image of Morse polynomials in
the highly composite case 
follows from the following 
criteria
(see \cite{crtpoisson}, Theorem~1):
%
{\em 
Assume that  there exists $\epsilon>0$ such that for each  integer
$k \geq 2$, 
\begin{equation}
  \label{eq:0}
N_k((h_1, h_2, \ldots, h_{k-1}),p) =  
\frac{p}{s_p^k} \left(1+ O_k( (1-s_p^{-1})p^{-\epsilon})\right)
\end{equation}
provided that  $0, h_1, h_2, \ldots, h_{k-1}$ are distinct mod $p$.
If $s_p = p^{o(1)}$ for all primes $p$, then the spacings modulo $q$
become Poisson distributed  as
$s_q$,  the mean spacing modulo $q$, tends to infinity.
}

What about non-Morse polynomials?  Rather surprisingly,
it turns out that 
~(\ref{eq:Nk-bound}) 
does not hold for all polynomials\footnote{In particular, the spacing
distribution for the image of such polynomials is {\em not} consistent
with the coin flip model!  (That is, independent coin flips where the
  probability of heads is given by $|\Sp|/p$.)}.
%
For example, in \cite{crtpoisson} it was shown that for $f(x) =x^4 -
2x^2$,  
$$
N_2(h,p) =
\begin{cases}
2/3 \cdot \frac{p}{s_p^2} + O(\sqrt{p}) 
&\text{ if $h \equiv \pm 1 \mod p$,  $p \equiv 1 \mod 4$}
\\
4/3 \cdot \frac{p}{s_p^2} + O(\sqrt{p}) 
&\text{ if $h \equiv \pm 1 \mod p$,  $p \equiv 3 \mod 4$}
\\
\frac{p}{s_p^2} + O(\sqrt{p}) 
&\text{ if $h \not \equiv \pm 1,0 \mod p$}
\end{cases}
$$
%
%

Hence the assumptions in (\ref{eq:0}) are violated.
However, we can prove that (\ref{eq:morse-mod-p-asymptotic})  holds for 
most values of $(h_0, h_1, \ldots, h_{k-1})$:
\begin{thm}
\label{thm:poly-correlation}
Let $p$ 
be a prime and let 
\begin{equation}
  \label{eq:R-p-definition}
R_p := \{ f(\xi) : f'(\xi) = 0, \xi \in \kbar \}.
\end{equation}
be the set of critical values modulo $p$.
If the sets $R_p, R_p-
h_1, R_p-h_2, \ldots, R_p-h_{k-1}$ are pairwise disjoint\footnote{In
  the case $f(x)=x^2$ this condition is equivalent to $0, h_1, \ldots,
  h_{k-1}$ being distinct modulo $p$.  However, for general 
  polynomials (including the case of Morse polynomials), the two conditions are
  {\em not} equivalent.}, then 
\begin{equation}
\label{eq:main-nk-asymptotic}
N_k((h_1, h_2, \ldots, h_{k-1}), p)
=
p/s_p^k + O_{f,k}(\sqrt{p})
\end{equation}
\end{thm}
In other words, the analogue of Davenport's result holds for all but
$O(p^{k-2})$ elements in $(\Z/p\Z)^{k-1}$.
Allowing for overlap between two translates of the set of critical values,
we also have the following weaker upper bound on $N_k((h_1, h_2, \ldots,
h_{k-1}), p)$:
\begin{prop}
\label{prop:poly-correlation-bound}
Let $p$ 
be a prime.
There exists a constant $C_0<1$, only depending on $f$, with
the following property:
if the sets  
$$
(R_p \cup R_p-h_1), R_p-h_2, \ldots, R_p-h_{k-1}
$$
are pairwise disjoint and $h_1 \not \equiv 0 \mod p$, then 
$$
N_k((h_1, h_2, \ldots, h_{k-1}), p)
\leq
\frac{C_0}{s_p^{k-1}} \cdot  p + O_{f,k}(\sqrt{p})
$$
\end{prop}
%
It turns out that these two results are enough to 
obtain Poisson spacings in the highly composite case.
However, rather than studying
the spacings 
directly, we proceed by determining the {\em $k$-level correlation
  functions}.
Let 
$$
\Sq
:=
\{ t \in \Z/q\Z : t = f(x) \text{ for some $x \in \Z/q\Z$} \}
$$ 
be the  image of $f$ modulo $q$,  let 
$$
s_q := q/|\Sq|
$$
be the mean spacing modulo $q$, and
given $\vh = (h_1, h_2, \ldots, h_{k-1}) \in \Z^{k-1}$, put
$$
N_k(\vh, q) := 
|\{ t \in \Sq : t+h_1, t+h_2, \ldots, t+h_{k-1} \in \Sq \}|
$$
For $X \subset \R^{k-1}$, the {\em $k$-level correlation function} is
then given 
by 
$$
R_k(X, q) := \frac{1}{|\Sq|} 
\sum_{\vh \in s_q X \cap \Z^{k-1}}
N_k(\vh,q)
$$
The main result of this paper
is then the following:
\begin{thm}
\label{thm:correlations}
Let $q$ be square free, $k \geq 2$ an integer, and let $X \subset
\R^{k-1}$ be a convex set with the property that $(x_0, x_1, \ldots
x_{k-1}) \in X$ implies that $x_i \neq x_j$ if $i \neq j$.  Then the
$k$-level correlation function of the image of $f$ modulo $q$
satisfies
$$
R_k(X, q) = \vol(X) + 
O_{f,k}\left( s_q^{-1/2+o(1)} + C_0^{\omega(q)(1-o(1))} \right)
$$
as $s_q \to \infty$, where $C_0<1$ is the constant given in
Proposition~\ref{prop:poly-correlation-bound}.
\end{thm}

Using a standard inclusion-exclusion argument (see \cite{squares1},
appendix A for details), this 
implies that the spacing statistics are
Poissonian.  In particular we have the following:
\begin{thm}
  For $q$ square free, the limiting (normalized)  spacing
  distribution\footnote{
By normalized spacings we mean the following: with
$0 \leq x_1<x_2< \cdots < x_{|\Sq|} <q$ 
being integer representatives of the
image of $f$ modulo $q$,  the spacings between consecutive elements
are defined to be $\Delta_i = x_{i+1} - x_i$ for $1 \leq i
< |\Sq|$,  and $\Delta_{|\Sq|} = x_1 - x_{|\Sq|} + q$.
The normalized spacings are then given by $\widetilde{\Delta_i} :=
\Delta_i/s_q$. 
} of
  the image of $f$ modulo $q$ is given by $P(t) = \exp(-t)$ as $s_q
  \to \infty$. 
  Moreover, for any integer $k \geq 2$, the limiting joint
  distribution of $k$ consecutive spacings is a product $\prod_{i=1}^k
  \exp( -t_i)$ of $k$ independent exponential variables.
\end{thm}

\subsection{Some remarks on the mean spacing}
\label{subsec:mean-spacing}
%
We note that the only way for which $s_p=1$ for all primes $p$ is if
$f(x)$ is of degree one.  However, there are nonlinear polynomials $f$
such that 
$s_p=1$ for infinitely many primes.  For example, if $f(x) = x^3$ and we take
$q$ to be a 
product of primes  
$p \equiv 2 \mod 3$, then $s_p=1$ for all $p|q$, and $s_q=\prod_{p|q}s_p = 1$
clearly does not tend to infinity.  
On the other hand, if $\deg(f)>1$, there is always a positive density 
set of primes $p$ such that $s_p>1$.  
Moreover,
if $f$ is not a permutation polynomial\footnote{$f$ is said to be a
  permutation polynomial 
modulo $p$ if $|\Sp| = p$.} modulo $p$,
Wan has shown
\cite{Wan} that 
\begin{equation}
\label{eq:wan}
|\Sp| \leq p -\frac{p-1}{\deg(f)}.  
\end{equation}
Thus, for primes $p$ such that $s_p>1$, $s_p$ is in fact uniformly
bounded away 
from $1$.  

It is also worth noting that  Birch and
Swinnerton-Dyer have shown \cite{bsd-valueset}
that for $f$ Morse, $|\Sp| = c_f \cdot p +
O_f(\sqrt{p})$ where $c_f<1$ only depends on the degree of $f$, hence
$s_p = 1/c_f + O(p^{-1/2})$ for all $p$, and thus $s_q \to \infty$ as
$\omega(q)\to \infty$.

\subsection{Related results}
There are only a few other cases for which Poisson spacings have been
proven.  Notable examples are Hooley's result \cite{Hooley2,Hooley3}
on invertible elements modulo $q$ under the assumption that the
average gap $s_q=q/\phi(q)$ tends to infinity, and the work by Cobeli
and Zaharescu \cite{cobeli-zaharescu} on spacings between primitive
roots modulo $p$, again under the assumption that the average gap $s_p
= (p-1)/\phi(p-1)$ tends to infinity.  Recently,
Cobeli,V{\^a}j{\^a}itu, and Zaharescu
\cite{cobeli-vaitu-zaharescu-inverses-mod-q} extended Hooley's results
and showed that subsets of the form $\{ x \mod q: x \in I_q, x^{-1}
\in J_q\}$ have limiting Poisson spacings if the intervals $I_q,J_q$ have large
lengths (more precisely, that $|I_q| \in [q^{1-(2/9(\log\log
  q)^{1/2})},q]$, and $|J_q| \in [q^{1-1/(\log \log q)^2}, q]$) as $q$
tends to infinity along a subsequence of integers such that $q/\phi(q) \to
\infty$.

\subsection{Acknowledgements}
The author would like to thank Juliusz Brzezi\'nski, Andrew Granville,
Moshe Jarden, Zeév Rudnick, and Thomas J. Tucker for  helpful discussions.



\section{Proof of Theorem~\ref{thm:poly-correlation}}
\label{sec:proof-of-thm-poly-correlation}

Given a polynomial $f \in \fp[x]$ and $k$ distinct elements $h_0=0,h_1,
h_2, \ldots, h_{k-1} \in \fp$, we wish to count the number of
$t \in \fp$ for which  there exists $x_0, x_1, \ldots x_{k-1}\in \fp$
such that
$$
f(x_0) = t+h_0, f(x_1) = t+h_1, f(x_2) = t+h_2, \ldots, f(x_{k-1}) =
t+h_{k-1}
$$
In order to study this, 
put 
$$\vh := (h_1, h_2, \ldots, h_{k-1})$$ 
and 
let $X_{k,\vh}$ be the affine curve defined by
$$
X_{k,\vh} := \{f(x_0) = t, \ f(x_1) = t+h_1, \ldots, f(x_{k-1}) =
t+h_{k-1}    \},
$$
and let $\fp[X_{k,\vh}]$ be the coordinate ring of $X_{k,\vh}$.  We
then have 
\begin{multline}
N_k((h_1, h_2, \ldots, h_{k-1}), \Sp)
\\
=
|\{\m \in \fp[t] : \text{$\M|\m$ for some degree one prime $\M \in
  \fp[X_{k,\vh}]$ } \}|
\end{multline}
In order to estimate the size of this set, we will use the Chebotarev
density theorem, made effective via the Riemann hypothesis for curves,
for the Galois closure of $\fp[X_{k,\vh}]$.  Thus, let 
$Y_{k,\vh}$ be the curve whose function field $\fp(Y_{k,\vh})$
corresponds to the Galois 
closure of the extension $\fp(X_{k,\vh})/\fp(t)$.  

We begin with the case $k=1$.  Given $h \in \kk$,
define a polynomial  $F_h \in \kk[x,t]$ by
$$
F_h(x,t) := f(x) - (t+h).
$$
Since the $t$-degree of $F_h$ is one,  $F_h$ is
irreducible, and thus $K_h = \kk[x,t]/F_h(x,t)$ is a field.
Let $L_h$ be the Galois closure of $K_h$, and let $G_h =
\gal(L_h/\kk(t))$ be the Galois group of the field extension
$L_h/\kk(t)$.   
By allowing for worse constants in the error terms, we may assume that
$p>n$, so that all
field extensions are separable, and no wild ramification can occur.

The following Lemma shows that $G_h$ and $L_h
\cap \overline{\kk}$ are independent of $h$.
\begin{lem}
\label{lem:same-constant-field}
Let $h \in \fp$.  Then
$
G_h \cong G_0
$
and $L_h \cap \overline{\kk} = L_0 \cap \overline{\kk}$.
\end{lem}
\begin{proof}
Define a $\kk$-linear automorphism $\sigma : \kk[t] \to \kk[t]$ by
$\sigma(t) = t + h$.
Since $\sigma(F_0) = F_h$ we may extend $\sigma$ to an
isomorphism $\sigma' : L_0 \to L_h$.  Moreover, given $\tau \in
G_0$, $\sigma' \tau (\sigma')^{-1} \in G_h$, the map $\tau \to
\sigma' \tau (\sigma')^{-1}$ gives an isomorphism between $G_0$ and
$G_h$.

Let $l_0 = L_0 \cap \overline{\kk}$ and let $l_h = L_h \cap
\overline{\kk}$. Since $l_0/\kk$ is normal, $l_0 = \sigma'(l_0) \subset
K_h \cap \kbar = l_h$, and the same argument for
$(\sigma')^{-1}$ gives that $l_h \subset l_0$, hence $l_h=l_0$.
\end{proof}

Thus  $$l := L_0 \cap \overline{\kk}$$  is the field
of constants for  $L_h$ for any $h \in \kk$.  Arguing as in the proof of
Lemma~\ref{lem:same-constant-field} we obtain:
\begin{lem}
Let $H_h := \gal(L_h/l(t))$.
Then
$
H_h \cong H_0.
$
\end{lem}

Our next goal is to obtain a criterion for linear disjointness for the
field extensions $L_h/l(t)$ as $h$ varies.
\begin{lem}
\label{lem:disjoint-ramification}
Let $E_1,E_2$ be finite extensions of $\kk(t)$, both having the same
constant field $l$, and degree smaller than $p$.
If $E_1/l(t)$ and $E_2/l(t)$ have disjoint finite ramification, then
$E_1 \cap E_2 = l(t)$.
\end{lem}
\begin{proof}
Let $E = E_1 \cap E_2$.  By the assumption, $E/l(t)$ can only ramify
at infinity.  Moreover, the ramification must be tame.
With $g_E$ denoting the genus of $E$, the Riemann-Hurwitz genus 
formula now gives
$$
-2 \leq 2(g_E-1) = [E:l(t)] 2( 0 -  1) +\sum_{\mathfrak{P}|\infty }
(e(\mathfrak{P}/\infty)-1) \deg(\mathfrak{P})
$$
$$
= - 2 [E:l(t)] + [E:l(t)] - \sum_{\mathfrak{P}|\infty }
\deg(\mathfrak{P})
< - [E:l(t)]
$$
and thus $ [E:l(t)] < 2$.
\end{proof}
We now easily obtain the desired criteria for linear disjointedness. 
\begin{prop}
If the sets $R_p, R_p-h_1$,
$R_p-h_2$, \ldots, $R_p - h_j$ are pairwise disjoint, then the field
extensions $L_0/l(t), L_{h_1}/l(t), \ldots,
L_{h_j}/l(t)$ are linearly disjoint.
\end{prop}
\begin{proof}
  Since $L_h$ is the Galois closure of $K_h$, both extensions,
  relative $\kk(t)$, ramify over the same primes.
The assumption of pairwise disjointness of $R_p, R_p-h_1,
\ldots, R_p- h_j$ means that there is no common finite
ramification among the fields $L_0, L_{h_1}, \ldots L_{h_j}$,
hence any intersection of compositums of the fields must, by
Lemma~\ref{lem:disjoint-ramification} and
Lemma~\ref{lem:same-constant-field}, equal $l(t)$
\end{proof}

If $G = \gal(E/\kk(t))$ is the Galois group of an extension
$E/\kk(t)$ with constant field $l$, define (following Cohen
\cite{cohen-poly-value-sets-one,cohen-poly-value-sets-two}) 
$$
G^*
:= \{ \sigma \in G :   \sigma |_{l(t)} =
\frob(l(t)/\kk(t)) \}
$$
where $\frob(l(t)/\kk(t))$ is the canonical generator of 
$\gal(l(t)/\kk(t))$ given by $x \to x^p$.
%

Let $L^k = \fp(Y_{k,\vh})$ be the
compositum of the fields $ L_{h_0}, L_{h_1},
\ldots, L_{h_{k-1}}$.
For $k \geq 2$, define a conjugacy class $\fix_{k,\vh} \subset
\gal(L^k/\kk(t))^*$ by 
\begin{multline*}
\fix_{k,\vh} := \{ \sigma \in \gal(L^k/k(t))^* : \\
\text{$\sigma$ fixes at least one root of $F_{h_i}$ for $i=0,
  1,\ldots, k-1$}
\}
\end{multline*}
For $k=1$ we define (note that there is no dependence on $\vh$)
a conjugacy class $\fix_{1} \subset \gal(L^1/\kk(t))^*$
by
$$
\fix_{1} := \{ \sigma \in \gal(L^1/k(t))^* : \\
\text{$\sigma$ fixes at least one root of $f(x)-t$} \}
$$

Then, taking into account $O_{k,f}(1)$  ramified primes, we have
\begin{multline}
N_k(\vh, p)
= \\ =
|\{\m \in \fp[t] : \text{$\deg(\m)=1$, $\exists \M|\m$, 
  $\M \subset \fp[Y_{k,\vh}]$,
$\frob(\M|\m) \in \fix_{k,\vh} $}  \}|
+O_{k,f}(1)
\end{multline}
where $\frob(\M|\m) \in \gal(L^k/\fp(t))$ denotes the Frobenius
automorphism. 
Applying the Chebotarev density theorem (e.g., see
\cite{fried-jarden-field-arithmetic}, Proposition~5.16), we obtain
$$
N_k(\vh, p)
=
\frac{|\fix_{k,\vh}|}{|\gal(L^k/l(t))|} \cdot p + O_{k,f}(\sqrt{p})
$$

Our next goal is to determine $|\fix_{k,\vh}|/|\gal(L^k/l(t))|$.
\begin{lem}
\label{lem:product-of-galois-group}
Given $k \geq 2$, define
$$
C_k(\vh,p)  := \frac{|\fix_{k,\vh}|}{|\gal(L^k/l(t))|},
$$
and 
$$
C_1(p) := \frac{|\fix_{1}|}{|\gal(L^1/l(t))|}.
$$
Assume that $R_p, R_p-h_1, \ldots, R_p-h_{k-1}$ are pairwise disjoint.
Then $C_k(\vh,p)= C_1(p)^k$ where $C_1(p)  = 1/s_p+ O_{f}(p^{-1/2})$.
\end{lem}

\begin{proof}
For simplicity\rmkd{Add $k>2$ case?}, we consider only the case $k=2$,
and for ease of notation, let $\vh = (h_1) = (h)$.

The action of $\gal(L^2/\kk(t))$ on the roots of $F_0$ and $F_h$
allows us to identify $\gal(L^2/\kk(t))$ and $\gal(L^2/l(t))$ with subgroups
of $S_n \times S_n$.  Moreover, since $L_0$ and $L_h$ are
linearly disjoint over $l(t)$ and have isomorphic Galois groups, we
may identify $\gal(L^2/l(t)) \cong H_0 \times H_h$ with a subgroup
of $S_n \times S_n$ in such a way that
$$
H_0 \cong
H' \times 1
\subset S_n \times 1
\subset S_n \times S_n
$$
and
$$
H_h \cong
1 \times H'
\subset 1 \times S_n
\subset S_n \times S_n
$$
where $H' \cong H_0 \cong H_h$ and $H'$ is a subgroup of $S_n$.

Define a $\kk$-linear map $\tau : \kk(t) \to \kk(t)$ by $\tau(t) = t+
h$, and extend it to a map from $L_0$ to $L_h$.
Given $\mu_1 \in G_0^*$, let $\mu_2 = \tau \mu_1 \tau^{-1}$.
Clearly $\mu_2 \in G_h$, and since $\gal(l(t)/\kk(t)) \cong
\gal(l/\kk)$ is abelian, 
$\mu_1|_{l(t)} = \mu_2|_{l(t)}$ and hence $\mu_2 \in G_h^*$.
Let us consider the possible extensions of $\mu_1,\mu_2$ to $L^2$.
After making a fixed, but arbitrary choice, of extensions
$\tilde{\mu_1}, \tilde{\mu_2}$ we
find that all pairs extensions are of the form $(\delta \mu_1, \gamma
\mu_2)$ where $\delta \in H_h$ and $\gamma \in H_0$.  Now, for
any such pair of extensions, we have
$$
\delta \mu_1 (\gamma \mu_2)^{-1} =
\delta \mu_1 \mu_2^{-1}  \gamma^{-1} \in \gal(L^2/l(t))
$$
But since $\gal(L^2/l(t)) \cong H_0 \times H_h$ we may choose
$\gamma$ and $\delta$ in such a way
that $\delta \tilde{\mu_1} \tilde{\mu_2}^{-1} \gamma^{-1} =1$.  In
other words, it is possible to choose $\tilde{\mu_1}, \tilde{\mu_2}$
so that $\tilde{\mu_1} = \tilde{\mu_2}$.

Thus, there is an extension of $\mu \in G_0^*$ to an element
$\tilde{\mu}$ of
$\gal(L^2/\kk(t))^*$
in such a way that $\tilde{\mu}$ embeds diagonally when regarded as an element
of $S_n \times S_n$, i.e., there exists $\sigma \in S_n$ such that
$\tilde{\mu}$ corresponds to
$$
(\sigma, \sigma) \in S_n \times S_n
$$

Now, all elements of $\gal(L^2/\kk(t))^*$, regarded as elements of
$S_n \times S_n$, 
must be of the form
$$
(\delta \sigma, \gamma \sigma) \in S_n \times S_n
$$
where $\delta,\gamma \in H'$.
In particular, if we let $H'' \subset H'$ be the set of elements
$\delta$ such that
$\delta \sigma$ has at least one fix point, we find that
$$
C_2(\vh,p) =
\frac{|H''|^2}{|\gal(L^2/l(t))|} =
\frac{|H''|^2}{|\gal(L^1/l(t))|^2} =
C_1(p)^2
$$
since $\gal(L^2/l(t)) \cong H_0 \times H_h$ and
$H_h \cong H_0 = \gal(L^1/l(t))$.

Finally, we note  that
$$
|\Omega_p| = p/s_p = |\{t \in \kk \text{ for which there exists $x
\in \kk$ such that $f(x)=t $}\}|
$$
$$
= C_1(p) \cdot  p + O_f(\sqrt{p})
$$
and thus $C_1(p)  = 1/s_p +O_f(\sqrt{p})$.
\end{proof}

\section{Proof of Proposition~\ref{prop:poly-correlation-bound}}
\label{sec:proof-prop-poly-correlation-bound}
We will begin by giving a proof for the case $k=2$, and then show how
the general case can be reduced to this case.  We will be using the
same notation as in the proof of Theorem~\ref{thm:poly-correlation},
and, by allowing worse constants in the error terms as before, we may
assume that $p>\deg(f)$.

\subsection{The case $k=2$}
We start by showing that the field extensions $K_0, K_h$ are disjoint
if $h \neq 0$.
\begin{lem}\label{l:absolutely-irreducible}
Let $f \in \kk[x]$ be a
polynomial of degree smaller  than $p$.  Then the affine curve
defined by 
$$
\{x,y: f(x) - (f(y)+h) = 0 \}
$$
is absolutely irreducible if $h \not \equiv  0 \mod p$.
\end{lem}
\begin{proof}
Let $x,y$ be roots of $f(x) = t$ and $f(y) = t+h$ where $t$ is
transcendental over $\kk$.  
If $\overline{\kk}(x)$ and $\overline{\kk}(y)$ are not linearly disjoint over
$\overline{\kk}(t)$ then, by Lüroth's theorem, 
$\overline{\kk}(x) \cap \overline{\kk}(y) = \overline{\kk}(u)$  for
some $u \not 
\in \overline{\kk}(t)$.  
Hence there exists non-constant rational functions $g_1$ and $g_2$
such that $u = g_1(x) = g_2(y)$, and a rational function $q$, of
degree less than $p$, such that $q(u)= t$.  However, since $t = f(x) =
q(u)=q(g_1(x))$ and $f$ is a polynomial, $q$ and $g_1$ must be of a
special form: either $q$ and $g_1$ are both polynomials, or $g_1(x)=
c_1 + c_2/b(x)$ where $c_1,c_2$ are constants, $b(x)$ is a polynomial,
and $q(u) = \sum_{i=0}^l a_i/(u-c_1)^i$.  In the latter case, we can
replace $u$ by $\tilde{u}= c_2/(c_1-u)$, and hence we may assume that
$q$ and $g_1$ are in fact both polynomials.  Similarly, since $q(g_2(y)) =
t = f(y)-h$, we may assume that $g_2$ is a polynomial as well.

%
Now, since $t$ is trancendental, so is $y$ and therefore
$q(g_2(y)) = t = f(y) - h $ implies that $q(g_2(x)) = f(x) - h $.
Thus 
$$
q(g_1(x)) - q(g_2(x)) = h
$$
and hence $g_1(x)-g_2(x)$ must divide $h$, which can only happen if 
$g_1(x)= g_2(x)+C$ for some  constant $C \neq 0$.  Thus
$$
q(g_2(x)+C) - q(g_2(x)) = h
$$
and hence
$$
q'(g_2(x)+C)g_2'(x) - q'(g_2(x))g_2'(x) = 0
$$
which, since $g_2$ is non-constant, implies that
$$
q'(g_2(x)+C) = q'(g_2(x))
$$
Therefore, if $g_2(\alpha) = \beta$ where $q'(\beta)=0$ we find that
$q'(\beta + C) = q'(\beta ) = 0$, and more generally, that 
$q'(\beta + lC) = 0$ for $l = 0, 1, \ldots p-1$, which is impossible
since the degree of $q$ is smaller than $p$.

Thus, the two fields  $\kbar(x,t)/(f(x)-t)$ and
$\kbar(y,t)/(f(y)-t-h)$ are linearly disjoint 
over $\kbar(t)$
and hence $f(x)-(f(y)+h)$, when
regarded as a polynomial over $\kbar(y)$, is irreducible.
\end{proof}

We are now ready to give a proof for
Proposition~\ref{prop:poly-correlation-bound} 
in the case $k=2$.
\begin{lem}
\label{lem:poly-correlation-bound-special-case}
There exists\rmkd{change notation to avoid $C_0$?} $\cnot < 1$, only
depending on $f$, with the following property:
for all sufficiently large $p$ for which 
$f$ is not a permutation polynomial modulo $p$,
$$
C_2((h),p) \leq \cnot /s_p
$$
if $h \not \equiv  0 \mod p$.

\end{lem}
\begin{proof}
For $f$ fixed there are only finitely many possibilities for
$\gal(L^2/\fp(t))$, hence $C_2((h),p) =
|\fix_{2,(h)}|/|\gal(L^2/l(t))|$ can only take finitely many
values.  Thus,  since
$C_2((h),p) \leq C_1(p) = 1/s_p + O_f(p^{-1/2})$
it is enough to show that 
$C_2((h),p) = C_1(p)$ can only happen for finitely many
primes $p$.

Given $a \in \fp$, let $M(a) = |\{x \in \fp : f(x) = a \}|$.  Then 
$$
|\{x,y \in \fp : f(x) = f(y) + h  \}|
=
\sum_{a \in \fp} M(a) M(a+h)
$$
On the other hand, by Lemma~\ref{l:absolutely-irreducible}, the curve
defined by $f(x) = 
f(y) + h$ is absolutely irreducible, and hence the 
Riemann hypothesis for curves gives that 
$$
|\{x,y \in \fp : f(x) = f(y) + h  \}|
=
p + O_f(\sqrt{p})
$$
We have
$$
|\{ a : M(a) > 0  \}|
=
|\{ a : M(a-h) > 0  \}|
 = |\text{Image}(f)|  = p/s_p
$$
Thus, if
\begin{multline*}
N_2(h,p) = |\{a \in \fp : M(a)>0, M(a-h)>0\}|= 
\\
C_2(h,p) \cdot p + O_f(\sqrt{p})
= C_1(p) \cdot p + O_f(\sqrt{p})
= \frac{1}{s_p} \cdot p + O_f(\sqrt{p})
\end{multline*}
then, since $ |\{ a : M(a-h) > 0  \}| = |\text{Image}(f)| = p/s_p $, we have
$$
|\{a \in \fp : M(a)=0, M(a-h)>0\}|
=  O_f(\sqrt{p})
$$
Therefore
\begin{multline*}
p+O_f(\sqrt{p}) = 
\sum_{a \in \fp} 
M(a)M(a-h) 
\\
\geq
\sum_{a \in \fp : M(a) = 1} 
M(a-h) 
+ 
2 \sum_{a \in \fp : M(a) > 1} 
M(a-h)
\\ =
\sum_{a \in \fp : M(a) > 0} 
M(a-h) 
+ 
\sum_{a \in \fp : M(a) > 1} 
M(a-h) 
\\ =
\sum_{a \in \fp }
M(a-h) 
+
\sum_{a \in \fp : M(a) > 1} 
M(a-h) 
-
\sum_{a \in \fp : M(a) = 0} 
M(a-h) 
\\ = 
p 
+
\sum_{a \in \fp : M(a) > 1} 
M(a-h) 
-
O_f(\sqrt{p})
\end{multline*}
and thus
$$
\sum_{a \in \fp : M(a) > 1} M(a-h) = O_f(\sqrt{p})
$$
Hence
$$
|\{ a \in \fp : M(a) > 1,  M(a-h)> 0\}|
= O_f(\sqrt{p})
$$
and we similarly obtain that
$$
|\{ a \in \fp : M(a) > 0,  M(a-h)> 1\}|
= O_f(\sqrt{p})
$$
But then
\begin{multline*}
p + O_f(\sqrt{p}) =
\sum_{a \in \fp} 
M(a)M(a-h) 
\\ =
|\{a \in \fp : M(a) = M(a-h) = 1\}|  +  O_f(\sqrt{p})
\end{multline*}
In other words, $M(a)=1$ for all but $O_f(\sqrt{p})$ elements, 
which, by Wan's result (see (\ref{eq:wan}),
section~\ref{subsec:mean-spacing}),  
can only happen if $f$ is bijection once  $p$ is sufficiently large.

\end{proof}

\subsection{The case $k>2$}
\label{sec:case-k2}
As usual, we use the convention that $h_0 = 0$.  
Arguing as in the proof of
Lemma~\ref{lem:disjoint-ramification}, we
find that the field extensions 
$$
\left(L_{h_0}L_{h_1}\right)/l(t), L_{h_2}/l(t), \ldots, L_{h_{k-2}}/l(t),L_{h_{k-1}}/l(t)
$$
are linearly disjoint since they have disjoint ramification.
Hence there is an isomorphism
\begin{multline*}
\gal \left(
L_{h_0} L_{h_1} \ldots \ldots L_{h_{k-1}}/l(t)
\right)
\\ \simeq
\gal\left(L_{h_0} L_{h_1}   /l(t)\right) \times
\gal\left( L_{h_2}/ l(t) \right) 
\times \ldots \times 
\gal\left( L_{h_{k-1}} / l(t) \right) 
\end{multline*}
Putting $\vh' = (h_0, h_1)$ and arguing as in
Lemma~\ref{lem:product-of-galois-group}, we find that  
$$
\frac{|\fix_{k,\vh}|}{|\gal(L^{k}/l(t))|}
=
\frac{|\fix_{2,\vh'}|}
{|\gal(L_{h_0}L_{h_1}/l(t))|}
\cdot \frac{1}{s_p^{k-2}}
= 
C_2(\vh',p) 
\cdot \frac{1}{s_p^{k-2}}.
$$
By Lemma~\ref{lem:poly-correlation-bound-special-case}, $C_2(\vh',p) \leq
C_0/s_p$ and the proof is complete.

\section{Proof of Theorem~\ref{thm:correlations}}

For $\vh \in \Z^{k-1}$ fixed, it follows immediately from the Chinese
Remainder Theorem that $N_k(\vh,q)$ is multiplicative in $q$.
The following Lemma shows that we
may assume that $q$ is a product of primes $p$ for which $f$ is not a
permutation polynomial modulo $p$, and hence that 
$s_p$ is uniformly bounded away from $1$ for all $p|q$.
\begin{lem} 
\label{lem:boring-permutations}
Given a square free integer $q$, write $q = q_1q_2$ where 
$$
q_1 = \prod_{\substack{p|q \\ |\Sp|< p}} p, \quad
q_2 = \prod_{\substack{p|q \\ |\Sp| = p}} p
$$
Then
$$
R_k(X, q) = R_k(X,q_1)
$$  
\end{lem}
\begin{proof}
If $p|q_2$ we have $s_p = p/|\Sp| = 1 $ and
$
N_k(\vh, p) =  p
$
for all $\vh \in \Z^{k-1}$.  
Thus $s_q = s_{q_1} \cdot s_{q_1} = s_{q_1}$, and since for $\vh$ fixed,
$N_k(\vh,q)$ 
is multiplicative, we find that $N_k(\vh,q) = N_k(\vh,q_1)\cdot q_2$.  
Thus
$$
R_k(X, q) = \frac{1}{|\Sq|} 
\sum_{\vh \in s_q X \cap \Z^{k-1}}
N_k(\vh,q)
=
\frac{q_2}{|\Sa||\Sb|} 
\sum_{\vh \in s_{q} X \cap \Z^{k-1}}
N_k(\vh,q_1)
$$
$$
=
\frac{1}{|\Sa|} 
\sum_{\vh \in s_{q_1} X \cap \Z^{k-1}}
N_k(\vh,q_1)
=
R_k(X,q_1)
$$

\end{proof}
We also note the following easy consequence of
Theorem~\ref{thm:poly-correlation}. 
\begin{lem}
\label{lem:weak-poly-correlation-bound}
Let $l$ be the largest integer
such that $R_p-h_{i_1}, R_p-h_{i_2},\ldots, R_p-h_{i_l}$ are pairwise
disjoint for some choice of indices $0 \leq i_1, i_2, \ldots, i_l \leq
k-1$ (with the usual convention that $h_0 = 0$).  Then
$$
N_k((h_1, h_2, \ldots, h_{k-1}), p)
\leq
p/s_p^{l}  + O_{f,k}(\sqrt{p})
$$
\end{lem}
\begin{proof}
If  $\{h_1', h_2', \ldots
h'_{l-1}\}$ is a subset of $\{h_1, h_2, \ldots, h_{k-1}\}$
then trivially
$$
N_k((h_1, h_2, \ldots, h_{k-1}), p) \leq
N_l((h_1', h_2', \ldots, h'_{l-1}), p)
$$
and the Lemma follows from Theorem~\ref{thm:poly-correlation}.
\end{proof}
%

\subsection{Some remarks on affine sets}

We will partition $\Z^{k-1}$ according to the size of the bounds on
$N_k(\vh, q) = \prod_{p|q} N_k(\vh, p)$ given by
Theorem~\ref{thm:poly-correlation} and
Proposition~\ref{prop:poly-correlation-bound}. 
In order to do this, we need to
introduce some notation:
By an {\em affine set} $L \subset \Z^{k-1}$ we mean an integer 
translate of a lattice $L' 
\subset \Z^{k-1}$.  We then define the rank, respectively discriminant, of
$L$ as the rank, respectively discriminant\footnote{By the
  discriminant of $L' \subset \Z^{k-1}$ we mean the index of $L'$ in
  $\Z^{k-1}$.}, of $L'$.  Similarly, we 
define $\codim(L)$ as $ k-1$ minus the rank of $L$.

Let $R$ be the set of critical values of $f$, i.e.,
$$
R := \{ f(\xi) : f'(\xi) = 0, \xi \in \overline{\Q} \}
$$
and recall that 
$
R_p = \{ f(\xi) : f'(\xi) = 0, \xi \in \kbar \}
$
is the set of critical values of $f$ modulo $p$.
Let 
$$
\sym := R - R = \{  \alpha - \beta : \alpha,\beta \in R  \},
$$
put 
$$\sym_\infty := \sym \cap \Z,$$ 
and let 
$$\sym_p := (R_p - R_p)  \cap \fp.$$ 
If $R_p + h_i \cap R_p+h_j \neq \emptyset$ then $h_i-h_j \in \sym_p$,
so the affine sets to be considered will be given by equations of the
form 
\begin{equation}
  \label{eq:affine-equality}
h_i - h_j = r, \ r \in \sym_\infty  
\end{equation}
or congruences of the form 
\begin{equation}
  \label{eq:affine-congruence}
h_i - h_j \equiv  r_p \mod p, \ r_p \in \sym_p  
\end{equation}
We note that the  bounds given by 
Theorem~\ref{thm:poly-correlation} and
Proposition~\ref{prop:poly-correlation-bound}
only depends on the congruence class of $\vh$, but we will
treat the case of equality separately since $N_k(\vh,p)$ will be large
{\em for all $p|q$} if $\vh$ satisfies an equation of the form
(\ref{eq:affine-equality}).

To ensure that the equations defining the affine sets are independent,
we will need the following notions:  Given $$E \subset \{(i,j) : 0 \leq
i<j \leq k-1\}$$ we may associate a 
graph $G(E)$ on the set of vertices $\{0, 1, \ldots, k-1\}$ by
regarding $E$ as the set of edges, i.e., two nodes $i,j$ are connected
by an edge if and only if $(i,j) \in E$.
Let 
$$
\AG := 
\{E \subset \{(i,j) : 0 \leq i<j \leq k-1\} : \text{ $G(E)$ is acyclic.} \}
$$
be the collection of  edge sets whose associated graphs are acyclic.

Given $E \in \AG$ and a map $\alpha : E \to \sym_\infty$, define an affine set
$$
L(E,\alpha) := \{ \vh \in \Z^{k-1} : h_i - h_j = \alpha((i,j)) \text{
  for all $(i,j) \in E$.} \}.
$$
(with the usual convention that $h_0=0$).  Note that $G(E)$ acyclic
implies that the equations defining $L(E,\alpha)$ are independent.
Further, given $E \in \AG$, let 
$$
\LL(E) := \{ L(E,\alpha) \text{ where $\alpha$ 
ranges over all maps $\alpha : E \to \sym_\infty$} \} 
$$
be the collection of affine sets defined by independent relations
between $h_i$ and $h_j$ for all $(i,j) \in E$.
We note that $\LL(\emptyset)$ contains exactly one element, namely the
full lattice $L(\emptyset, - ) =\Z^{k-1}$.
Moreover,  if $L \in \LL(E)$, then (since we assume that $E \in \AG$)
$\codim(L) = |E|$, and if $\vh \in L$, then
Proposition~\ref{prop:poly-correlation-bound} will, for {\em all} $p|q$, at 
best give the bound
$$
N_k(\vh,p ) \leq C_0 \frac{p}{s_p^{k-|E|}} + O_{f,k}(\sqrt{p}).
$$
(The bound will not hold if the components of $\vh$ satisfies
additional equations, i.e., if $\vh \in L'$ for some $L' \in \LL(E')$
such that $E' \supsetneq E$.)


Given $L(E, \alpha) \in \LL(E)$, let
$$
L^\times(E, \alpha) := 
\{ \vh \in L(E, \alpha) : \vh \not \in L(E',\alpha') 
\text{ for all $E' \supsetneq E$, $\alpha' : E' \to \sym_\infty$}\}
$$
In particular, if $\vh \in L^\times(E, \alpha)$, the components of
$\vh$ satisfy exactly $|E|$ independent 
  equations of the 
form $h_i-h_j = r_{ij}$ where $r_{ij} \in \sym_\infty$.

We also need to keep track of similar relations, modulo $p$, between
the components of $\vh$.  
Thus, given $E_p \in \AG$ and $\alpha_p : E_p \to \sym_p$, define an
affine set
$$
L_p(E_p,\alpha_p) := \{ \vh \in \Z^{k-1} : h_i - h_j \equiv \alpha_p((i,j))
\mod p \text{  for all $(i,j) \in E_p$} \}.
$$
We note that the rank of $L_p(E_p,\alpha_p)$
is $k-1$ and that the discriminant of $L_p(E_p,\alpha_p)$ is $p^{|E_p|}$,
and if $\vh \in L_p(E_p, \alpha_p)$, then  
Proposition~\ref{prop:poly-correlation-bound} will at
best give the bound
$$
N_k(\vh,p ) \leq C_0 \frac{p}{s_p^{k-|E_p|}} +
O_{f,k}(\sqrt{p}).
$$
Now, given $E \in \AG$, let 
$$
\LL_p(E) := \{ L_p(E_p,\alpha_p) : E_p \in \AG, \alpha_p : E_p \to \sym_p,
E_p \cap E = \emptyset, 
E_p \cup E \in \AG \}
$$
and for
$L_p \in \LL_p(E)$, let 
$$
L_p^\times := 
\{ \vh \in L_p : \vh \not \in L_p' 
\text{ for all $L_p' \in \LL_p(E_p')$, $E_p' \supsetneq E_p$}
\}
$$
If $\vh \in L^\times \cap L_p^\times$ for $L \in \LL(E)$ and $L_p =
L_p(E_p, \alpha_p) \in \LL_p(E)$, then $(h_0, h_1, \ldots,
h_{k-1})=\vh$ satisfies exactly $|E|$ independent equations of the
form $h_i-h_j = r_{ij}$ where $r_{ij} \in \sym_\infty$, and exactly
$|E_p|$ independent congruences of the $h_i-h_j \equiv r'_{ij} \mod p$
where $r'_{ij} \in \sym_p$, and furthermore,  there is no overlap
between the  equations and congruences.
The reason for keeping track of equalities and congruences separately
is that if $\vh \in L$ for $L \in \LL(E)$ and $|E|>0$, then the bounds
given on $N_k(\vh,p)$ given by
Proposition~\ref{prop:poly-correlation-bound} allows $N_k(\vh,p)$ to
deviate quite a bit from its mean value {\em for all $p|q$}.  On the
other hand, if we let $c$ be the product of primes $p|q$ for which the
bounds are bad because of congruence conditions, rather than
equalities, then we can bound the size of $c$ (see
Lemma~\ref{lem:c-must-be-bounded}).  We can now partition $\Z^{k-1}$
according to the size of the bounds on $N_k(h,p)$ given by
Theorem~\ref{thm:poly-correlation} and
Proposition~\ref{prop:poly-correlation-bound}:

\begin{lem}
\label{lem:bound-for-L-p}
Let $L = \lattice(E,\alpha)$,  $L_p = L_p(E_p, \alpha_p) \in
\LL_p(E)$, and assume that $\vh \in L^\times \cap L_p^\times$. 
If $|E|+|E_p| = 0$, then 
$$
N_k(\vh,p)
=
s_p^{-k} \cdot p + O_{k,f}(p^{1/2}),
$$
whereas if  $k>|E|+|E_p| > 0$, then 
$$
N_k(\vh,p)
\leq 
C_0 \cdot s_p^{|E|+|E_p|-k} \cdot p + O_{k,f}(p^{1/2}).
$$
where $\cnot<1$ is as in Proposition~\ref{prop:poly-correlation-bound}.

\end{lem}
\begin{proof}
The first assertion follows immediately from
Theorem~\ref{thm:poly-correlation} since $R_p + h_i \cap R_p + h_j
\neq \emptyset$ implies that $h_i-h_j \in \sym_p$.

For the second assertion, we argue as follows: 
Since $\vh = (h_1,h_2, \ldots, h_{k-1}) \in \lstar \cap \lstar_p$
there are indices 
$i_1, i_2, \ldots, i_{k-|E|-|E_p|}$ such that  $h_{i_1} \neq h_{i_2}$ and 
$$
(R_p - h_{i_1} \cup  R_p- h_{i_2}), R_p - h_{i_3}, \ldots, R_p - h_{i_{k-|E|-|E_p|}}
$$
are pairwise disjoint.  Putting 
$$
\vh' = (h_{i_2}-h_{i_1}, h_{i_3} -h_{i_1} , \ldots,
h_{i_{k-|E|-|E_p|}}-h_{i_1} ),
$$ 
the result follows from the bound
for $N_k(\vh',p)$ given by 
Proposition~\ref{prop:poly-correlation-bound}.
\end{proof}

However, partitioning $\Z^{k-1}$ according to the size of $N_k(\vh,p)$
for individual prime factors $p|q$ 
is not quite enough; we need to partition 
$\Z^{k-1}$ according to the size of $N_k(\vh,q) = \prod_{p|q}
N_k(\vh,p)$.  Thus, let  
$$
\LL_c(E) := 
\{ L \cap (\cap_{p|c} L_p) : 
L \in \LL(E),  
\ \forall p|c\ L_p \in \LL_p(E) \setminus L_p(\emptyset, -) \}
$$
(where $L_p( \emptyset, -  ) \in \LL_p(E)$ is the maximal lattice, i.e., 
$L_p( \emptyset, -  ) = \Z^{k-1}$) and given 
$$
L_c = L \cap (\cap_{p|c} L_p) \in \LL_c(E)
$$
let 
$$
L_c^\times :=
L^\times \cap (\cap_{p|c} L_p^\times) 
\cap (\cap_{p| \frac{q}{c}} L_p^\times( \emptyset, -  ) )
$$
We can now partition $\Z^{k-1}$ into subsets $L_c^\times$, where $L_c \in
\LL_c(E)$, $E \in \AG$, and $c|q$.
Moreover, as an
immediate consequence of the definitions and
Lemma~\ref{lem:bound-for-L-p}, we  
obtain the following:
\begin{lem}
\label{lem:bound-for-L-c}
Assume that $L_c = L \cap (\cap_{p|c} L_p(E_p, \alpha_p)) \in \LL_c(E)$ and
that $\vh \in L_c^\times$. 
If $p \nmid c$, then 
$$
N_k(\vh,p)
=
s_p^{-k} \cdot p + O_{k,f}(p^{1/2}).
$$
If $p \mid c$, then 
$$
N_k(\vh,p)
\leq
C_0 \cdot s_p^{|E|+|E_p|-k} \cdot p + O_{k,f}(p^{1/2}).
$$
where $\cnot<1$ is as in Proposition~\ref{prop:poly-correlation-bound}.

\end{lem}


Using the previous Lemma we can now bound sums of the form $\sum_{\vh \in
s_q X \cap \lstar_c} N_k(\vh, q)$.
\begin{lem}
\label{lem:lots-of-lc-bounds}
If 
$$
\lattice_c = L \cap (\cap_{p|c}  \lattice_p(E_p, \alpha_p)  )  \in \LL_c(E),
$$
then 
$$
|\{ \vh \in s_q X \cap \lstar_c  \}|
\leq 
|\{ \vh \in s_q X \cap \lattice_c  \}|
\ll_{k,f,X}
\frac{s_q^{k-|E|-1}}{c} + 
s_q^{k-|E|-2} 
$$  
Moreover, if $\vh \in \lstar_c$, then 
$$
\frac{N_k(\vh, q)}{q/s_q} \ll 
\prod_{p|c }\left( \frac{s_p^{|E|+|E_p|}}{s_p^{k-1}} + O_{k,f}(p^{-1/2}) \right)
\cdot 
\prod_{p|\frac{q}{c}}
\left( \cnot \cdot \frac{ s_p^{|E|}}{s_p^{k-1}} + O_{k,f}(p^{-1/2})  \right)
$$
In particular,
\begin{multline}
  \label{eq:why-truncate-c}
\sum_{\vh \in s_q X \cap \lstar_c} 
\frac{N_k(\vh, q)}{q/s_q} 
\\ \ll 
s_c^{k-1} \cnot^{-\omega(c)}
(\frac{1}{s_q}  + \frac{1}{c}  )
\cdot 
\cnot^{\omega(q)}
\cdot 
\prod_{p|q}
\left(  1 +  O_{k,f}(p^{-1/2})  \right)
\end{multline}

\end{lem}
\begin{proof}
The first assertion  follows from the Lipschitz
principle\footnote{Actually, we have to be a little 
careful: if we embed $\lattice$ into
$\Z^{k-1-|E|}$ and apply the Lipschitz principle,
there is an implicit constant in the bound that will depend on
$\lattice$.  However, the estimate is
uniform since $\lattice$ only can be
choosen in $O_k(1)$ ways.} (e.g.,
see Lemma~16 in \cite{squares1}) since $\lattice_c$ is a translate of
a lattice with discriminant (relative $L$) divisible by $c$.  The second
assertion follows from 
Lemma~\ref{lem:bound-for-L-c}.  Thus
\begin{multline*}
\sum_{\vh \in s_q X \cap \lstar_c} 
\frac{N_k(\vh, q)}{q/s_q} \ll 
\prod_{p|c }
\left( \frac{s_p^{|E_p|}}{p} + O_{k,f}(p^{-3/2}) \right)
\cdot 
\prod_{p|\frac{q}{c}}
\left(  \cnot +  O_{k,f}(p^{-1/2})  \right)
\\
+
\frac{1}{s_q}
\prod_{p|c }
\left( s_p^{|E_p|} + O_{k,f}(p^{-1/2}) \right)
\cdot 
\prod_{p|\frac{q}{c}}
\left(  \cnot +  O_{k,f}(p^{-1/2})  \right)
\end{multline*}
$$
\ll
 \cnot^{-\omega(c)}
( \frac{s_c^{k-1}}{c} + \frac{s_c^{k-1}}{s_q}   )
\cdot 
\cnot^{\omega(q)}
\cdot 
\prod_{p|q}
\left(  1 +  O_{k,f}(p^{-1/2})  \right)
$$

\end{proof}

Since the bound in (\ref{eq:why-truncate-c}) is not useful for large
$c$, we will also need the following:
\begin{lem}
\label{lem:c-must-be-bounded}
Let $d$ be the degree of the field extension $\Q(\sym)/\Q$.
If $\lattice_c \in \LL_c(E)$ for some $E \in \AG$
and  $s_q X \cap \lattice_c^\times  \neq \emptyset$ then 
\begin{equation*}
  \label{eq:small-c}
c \ll_{X,\sym} s_q^{d \binom{k}{2} |\sym|}.  
\end{equation*}
Moreover, there exist a constant $D$, only depending on $k$ and $f$, such 
that 
$$
|\LL_c(E)| 
\ll_{k,f} D^{\omega(c)}.
$$
\end{lem}
\begin{proof}
We first assume that all elements of $\sym$ are algebraic
integers.  
Let $B$ be the ring of integers in $\Q(\sym)$.  For each
prime $p|q$ chose a  prime ${\mathfrak P}_p \subset B$ lying
above $p$, so that we may regard any element in $\sym_p$ as  
the image of an element in $\sym$ under the reduction map $B \to
B/{\mathfrak P}_p$.

For $0 \leq i < j \leq k-1$, $r \in \sym$, and $\vh \in
\lattice_c^\times$, let 
$$
\gamma_{i,j,r}(\vh) = \prod_{p : h_i - h_j \equiv r \mod {\mathfrak P}_p} p
$$
Then $c$ divides 
$$
\prod_{\substack{ 0 \leq i < j \leq k-1 \\ 
r \in \sym : h_i-h_j \neq r }} \gamma_{i,j,r}(\vh)
$$
Since $h_i-h_j-r \equiv 0 \mod {\mathfrak P}_p$ for all $p$ dividing
$\gamma_{i,j,r}$, we find that $\gamma_{i,j,r}$ divides
$\norm_{\Q}^{\Q(\sym)}(h_i-h_j-r)$.  
Moreover,  if $\vh  \in s_q X$, then $|h_i -h_j| \ll_X s_q$, thus 
$$
\norm_{\Q}^{\Q(\sym)}(h_i-h_j-r) \ll_{f,X} s_q^d
$$
and hence 
$$
c \leq
\prod_{\substack{ 0 \leq i < j \leq k-1 \\ 
r \in \sym : h_i-h_j \neq r }} 
\norm_{\Q}^{\Q(\sym)}(h_i-h_j-r) 
\ll_{k,f,X}
s_q^{d |\sym| \binom{k}{2} }
$$
(Note that
$\norm_{\Q}^{\Q(\sym)}(h_i-h_j-r) \neq 0$ since $h_i-h_j-r \neq 0$).

In case $\sym$ contains elements that are not algebraic integers, we
can find an integer $m$, only depending on $\sym$, such that all
elements of $m \cdot \sym = \{ m\cdot r : r \in \sym \}$ are algebraic
integers, and apply the above argument to $m \cdot \sym$ and $m \vh$
(for primes $p$ not dividing $m$, but since $c$ is square free this
just makes the constant worse by a power of $(c,m) \leq  m$, which is
$O(1)$.) 

The second assertion follows upon noting that  there are
 $O_{k,f}(1)$ possible choices of $E_p$ and $\alpha_p$ for each $p|c$.

\end{proof}

\subsection{Conclusion}

We can now write $\Z^{k-1}$ as a disjoint union of sets
$\lstar$ where $\lattice$ ranges over all elements in $\cup_{E \in
  \AG} \LL(E)$, and hence $R_k(X,q)$ equals
\begin{multline}
\label{poly-eq:main-term}
\frac{1}{|\Sq|}
\sum_{\vh \in s_q X \cap \Z^{k-1}} 
N_k(\vh,q)
 =
\frac{1}{|\Sq|}
\sum_{E \in \AG}
\sum_{L \in \LL(E)}
\sum_{\vh \in s_q X \cap \lstar}
N_k(\vh,q)
\end{multline}
%
The term corresponding to $E=\emptyset$ in (\ref{poly-eq:main-term})
will give the main 
contribution
(note that if $ E = \emptyset$, then
$L = \lattice_\infty(E, -)= \Z^{k-1}$.)
Let 
$$
X' := \{\vh \in X : h_i - h_j \not \in \sym_\infty 
\text{ for $0\leq i < j \leq k-1$}  \}
$$
where we as usual use the convention that $h_0 = 0$.  Then
$$
s_q X \cap \lstar = s_q X' \cap \Z^{k-1}
$$
Note that  
$X'$ is just $\R^{k-1}$ with some 
hyperplanes removed, so if $X$ is convex, we can write $X'$ as a
finite union of convex sets.
We now rewrite (\ref{poly-eq:main-term}) as follows:
$$
\frac{1}{|\Sq|}
\sum_{\vh \in s_q X \cap \Z^{k-1}} 
N_k(\vh,q)
=
\sum_{\vh \in s_q X' \cap \Z^{k-1}} 
N_k(\vh,q)
+
\text{Error}_1
$$
where 
$$
\text{Error}_1 := 
\frac{1}{|\Sq|}
\sum_{E \in \AG, |E|>0}
\sum_{L \in \LL(E)}
\sum_{\vh \in s_q X \cap \lstar}
N_k(\vh,q)
$$
and the main term is given by
\begin{equation}
  \label{eq:main-term-sum}
\sum_{\vh \in s_q X' \cap \Z^{k-1}} 
N_k(\vh,q)
\end{equation}

We begin by showing that $\text{Error}_1 = o(1)$ as $\omega(q) \to
\infty$. 
\begin{lem}
As $\omega(q) \to \infty$, 
$$
\text{Error}_1 = 
\frac{1}{|\Sq|}
\sum_{\substack{E \in \AG \\|E|>0}}
\sum_{L \in \LL(E)}
\sum_{\vh \in s_q X \cap \lstar}
N_k(\vh,q)
\ll
\cnot^{\omega(q)(1-o(1))}.
$$
  
\end{lem}
\begin{proof}
Given $E \in \AG$ with $|E|>0$, we find that
\begin{multline}
  \label{poly-eq:6}
\frac{1}{|\Sq|}
\sum_{L \in \LL(E)}
\sum_{\vh \in s_q X \cap \lstar}
N_k(\vh,q)
\\ =
\frac{1}{q/s_q}
\sum_{c|q}
\sum_{\lattice_c \in \LL_c(E)}
\sum_{\vh \in s_q X \cap \lstar_c}
N_k(\vh,q)
\end{multline}
which, by Lemmas~\ref{lem:lots-of-lc-bounds} and
\ref{lem:c-must-be-bounded} is
\begin{equation}
  \label{poly-eq:5}
\ll
\cnot^{\omega(q)} \cdot 
\prod_{p|q}
\left(  1 +  O(p^{-1/2})  \right)
\sum_{\substack{c|q \\ c \ll s_q^\cboundexponent}}
D^{\omega(c)}
s_c^{k-1} \cnot^{-\omega(c)}
(\frac{1}{s_q}  + \frac{1}{c}  )
\end{equation}
Now,
$$
\sum_{\substack{c|q \\ c \ll s_q^\cboundexponent}}
D^{\omega(c)}
s_c^{k-1} \cnot^{-\omega(c)}
\frac{1}{c}  
\ll
\prod_{p|q}
\left( 1+O(1/p) \right)
$$
and, for any $\delta>0$, 
$$
\frac{1}{s_q}  
\sum_{\substack{c|q \\ c \ll s_q^\cboundexponent}}
D^{\omega(c)}
s_c^{k-1} \cnot^{-\omega(c)}
\ll
\frac{1}{s_q^{1-\delta \cboundexponent}}  
\sum_{c|q}
\frac{s_c^{k-1} \cnot^{-\omega(c)}}{c^\delta}
$$
$$
\ll
\frac{1}{s_q^{1-\delta \cboundexponent}}  
\prod_{p|q} 
\left( 1+O(1/p^\delta) \right)
\ll
\frac{1}{s_q^{1-\delta \cboundexponent-o(1)}}  
$$
Thus, taking $\delta = 1/(2\cboundexponent)$, we find that (\ref{poly-eq:5}) is 
$$
\ll 
\cnot^{\omega(q)} \cdot 
\prod_{p|q}
\left(  1 +  O(p^{-1/2})  \right)
\cdot 
\left(  
\frac{1}{s_q^{1/2-o(1)}}
+ 
\prod_{p|q}\left(  1 +  O(p^{-1})  \right)
  \right)
$$
$$
\ll
\cnot^{\omega(q)} \cdot 
\prod_{p|q}
\left(  1 +  O(p^{-1/2})  \right)
=
\cnot^{\omega(q)(1-o(1))}
$$
Since there are $O(1)$ possible choices of $L \in
\LL(E)$ for $E$ fixed, and
$E$ ranges over a finite number of subsets,  we find that
(\ref{poly-eq:6}) is $\cnot^{\omega(q)(1-o(1))}.$
  
\end{proof}

We procede by rewriting the main term  
in terms of a divisor sum.
%
For $p$ prime and $\vh \in \Z^{k-1}$, let
$$
\eps_k(\vh,p)=
\frac{s_p^{k-1} \cdot N_k(\vh,p)}{ |\Sp|} -1,
$$
so that we may write
$$
N_k(\vh,p)
= 
 \frac{|\Sp|}{s_p^{k-1}}(1 + \eps_k(\vh,p))   
$$
(recall that $s_p=p/|\Sp|$.) 
Further, for  $d>1$ a square free integer, put
$$
\eps_k(\vh,d) = \prod_{p|d} \eps_k(\vh,p) 
$$
and, to make $\eps_k$ multiplicative in the second parameter, set
$\eps_k(\vh,1) = 1$ for all $h$. 
%
%
Since $N_k(\vh,q)$ is multiplicative, we then have
\begin{equation}
N_k(\vh,q)
=
\prod_{p|q} 
\frac{1}{s_p^{k-1}}
|\Sp| 
\left(1 + \eps_k(\vh,p) \right)
=
\frac{|\Sq|}{s_q^{k-1}}
\sum_{d|q}
\eps_k(\vh,d)
\end{equation}

The following Lemma shows that the average of $\eps_k(\vh,d)$, over a
full set of residues modulo $d$, equals zero if $d>1$.
\begin{lem}
\label{lem:eps-average-equals-zero}
If $d>1$ then 
$$
\sum_{\vh \in (\Z/d\Z)^{k-1}} 
\eps_k(\vh,d)
=
0
$$
\end{lem}
\begin{proof}
Since $\eps_k(\vh,d)$ is multiplicative it is enough to show that
$$
\sum_{\vh \in (\Z/p\Z)^{k-1}} 
\eps_k(\vh,p)
= 0
$$
for $p$ prime, and because
$$
N_k(\vh,p)
= 
\frac{1}{s_p^{k-1}}
|\Sp|(1 + \eps_k(\vh,p))   
$$
it is enough to show that 
$$
\sum_{\vh \in (\Z/p\Z)^{k-1}} 
N_k(\vh,p)
= 
\frac{1}{s_p^{k-1}}|\Sp| p^{k-1}
=
|\Sp|^k
$$
But $\sum_{\vh \in (\Z/p\Z)^{k-1}} N_k(\vh,p)$ equals the number of
$k$-tuples of 
elements from $\Sp$, and hence  
$
\sum_{\vh \in (\Z/p\Z)^{k-1}} N_k(\vh,p)
=
|\Sp|^{k}.
$
\end{proof}

We will also need the following bound:
\begin{lem}
\label{lem:abs-eps-bound}
We have
$$
\sum_{\vh \in (\Z/d\Z)^{k-1}} 
|\eps_k(\vh,d)|
\ll 
d^{k-3/2+o(1)}
$$  
\end{lem}
\begin{proof}
Since the sum is multiplicative in $d$, it is enough to show that
$$
\sum_{\vh \in (\Z/p\Z)^{k-1}} 
|\eps_k(\vh,p)|
\ll 
p^{k-3/2}
$$
for $p$ prime.  By Theorem~\ref{thm:poly-correlation}, 
$|\eps_k(\vh,p)| \ll p^{-1/2}$
for all but $O(p^{k-2})$ residues modulo $p$, and for the remaining
residues we have $|\eps_k(\vh,p)| = O_{k,f}(1)$.  Thus 
$$
\sum_{\vh \in (\Z/p\Z)^{k-1}} 
|\eps_k(\vh,p)|
\ll 
p^{k-1} p^{-1/2} + p^{k-2}
\ll
p^{k-3/2}
$$
\end{proof}

We now find that the  main term (\ref{eq:main-term-sum}) equals
$$
\frac{1}{|\Sq|}  
\sum_{\vh \in s_q X' \cap \Z^{k-1}}
N_k(\vh,q)
=
\frac{1}{s_q^{k-1}}
\sum_{d|q}
\sum_{\vh \in s_q X' \cap \Z^{k-1}}
\eps_k(\vh,d)
$$
$$
= 
\frac{1}{s_q^{k-1}}
\sum_{\vh \in s_q X' \cap \Z^{k-1}}
1
+ \text{Error}_2
$$
where
\begin{equation*}
\text{Error}_2 :=
\frac{1}{s_q^{k-1}}
\sum_{\substack{ d|q \\ d>1}}
\sum_{\vh \in s_q X' \cap \Z^{k-1}}
\eps_k(\vh,d)
\end{equation*}
and the modified main term is
\begin{multline*}
  \frac{1}{s_q^{k-1}}
\sum_{\vh \in s_q X' \cap \Z^{k-1}}
1
=
\frac{1}{s_q^{k-1}}
\left( \vol(s_qX') + O(s_q^{k-2})  \right)
\\=
\vol(X) + O(1/s_q).
\end{multline*}

We conclude by  showing that  $\text{Error}_2 = o(1)$ as $s_q \to \infty$.
\begin{lem} As $s_q \to \infty$, we have
  \begin{equation}
    \label{eq:small-error}
\text{Error}_2 = 
\frac{1}{s_q^{k-1}}
\sum_{\substack{ d|q \\ d>1}}
\sum_{\vh \in s_q X' \cap \Z^{k-1}}
\eps_k(\vh,d)
\ll
s_q^{-1/2+o(1)}
  \end{equation}
\end{lem}
\begin{proof}

In order to show that $\text{Error}_2$ is small, we split the divisor
sum in two parts according to the size of $d$.

{\em Small $d$:}  We first consider $d \leq s_q^T$ where $T \in (0,1)$
is to be chosen later.
A point $\vh \in s_q X' \cap \Z^{k-1}$ is contained in
a unique  
cube $C_{\vh,d} \subset \R^{k-1}$ of the form
$$
C_{\vh,d} = \{(x_1, x_2, \ldots, x_{k-1}) : dt_i \leq x_i < d(t_i+1),
t_i \in \Z, \, i =1,2, \ldots, k-1  \}
$$
We say that $\vh \in s_q X' \cap \Z^{k-1}$ is a {\em $d$-interior} point
of $s_qX'$ if $C_{\vh,d} \subset 
s_q X'$, and if $C_{\vh,d}$ intersects the boundary of $s_qX'$, we say
that $h$ is a {\em $d$-boundary point} of $s_qX'$.

By Lemma~\ref{lem:eps-average-equals-zero}, the sum over the
$d$-interior points is zero, and hence
\begin{multline}
  \label{poly-eq:4}
\frac{1}{s_q^{k-1}}
\sum_{\substack{ d|q \\ 1< d \leq s_q^T}}
\sum_{\vh \in s_q X' \cap \Z^{k-1}}
\eps_k(\vh,d)
\\=
\frac{1}{s_q^{k-1}}
\sum_{\substack{ d|q \\ 1< d \leq s_q^T}}
\sum_{\substack{\vh \in s_q X' \cap \Z^{k-1}\\\text{$\vh$ is
      $d$-boundary point}}} 
\eps_k(\vh,d)
\end{multline}
Since $s_qX'$ is a union of convex sets, the number of cubes
$C_{\vh,d}$ intersecting the boundary of $s_qX'$ is $\ll (s_q/d)^{k-2}
$, and hence (\ref{poly-eq:4}) is
$$
\ll
\frac{1}{s_q^{k-1}}
\sum_{\substack{ d|q \\ 1< d \leq s_q^T}}
(s_q/d)^{k-2}
\sum_{\vh \in (\Z/d\Z)^{k-1}}
|\eps_k(\vh,d)|
$$
\begin{equation}
  \label{poly-eq:7}
=
\frac{1}{s_q}
\sum_{\substack{ d|q \\ 1< d \leq s_q^T}}
\frac{1}{d^{k-2}}
\sum_{\vh \in (\Z/d\Z)^{k-1}}
|\eps_k(\vh,d)|
\end{equation}
which by Lemma~\ref{lem:abs-eps-bound} is, for any $\alpha>1/2$, 
$$
\ll
\frac{1}{s_q}
\sum_{\substack{ d|q \\ 1< d \leq s_q^T}}
d^{1/2+o(1)}
\leq
s_q^{\alpha T-1}
\sum_{ d|q }
d^{1/2-\alpha+o(1)}
\ll
s_q^{\alpha T-1+o(1)}
$$
since
$$
\sum_{d|q} d^{-\epsilon} = \prod_{p|q}(1+p^{-\epsilon}) = s_q^{o(1)}
$$
if $\epsilon>0$ (recall that $s_p$ is
assumed to be uniformly bounded away from $1$
and $s_q = \prod_{p|q} s_p$.)

{\em Large $d$:}  We now consider
\begin{equation}
\label{poly-eq:2}
\frac{1}{s_q^{k-1}}
\sum_{\substack{ d|q \\ d>s_q^T}}
\sum_{\vh \in s_q X' \cap \Z^{k-1}}
\eps_k(\vh,d)
\end{equation}

Given $\vh$ and $d$, let $c$ be the largest divisor of $d$  such that
$\vh \in \lattice_c$ for some $\lattice_c \in \LL_c(L)$.
Then  
$$
\eps_k(\vh,d) 
\ll
\frac{s_c^{k-1}}{(d/c)^{1/2-o(1)}}
$$
by Lemma~\ref{lem:bound-for-L-c}.  Hence, for $E \in \AG$ fixed, 
$$
\sum_{L \in \LL(E)}
\sum_{\vh \in s_q X \cap \lstar}
\eps_k(\vh,d)
\ll
\sum_{c|d}
\sum_{\lattice_c \in \LL_c(E)}
\sum_{\vh \in s_q X \cap \lattice_c^\times}
|\eps_k(\vh,d)|
$$

$$
\ll
\sum_{c|d}
\frac{s_c^{k-1}}{(d/c)^{1/2-o(1)}}
\sum_{\lattice_c \in \LL_c(E)}
\sum_{\vh \in s_q X \cap \lattice_c^\times}
1
$$
which by Lemmas~\ref{lem:lots-of-lc-bounds} and \ref{lem:c-must-be-bounded}
is 
\begin{equation}
  \label{poly-eq:3}
\ll 
s_q^{k-1} \cdot 
d^{-1/2+o(1)}  \cdot 
\sum_{\substack{c|d \\ c \ll s^{\cboundexponent}}}
s_c^{k-1}c^{1/2-o(1)}
D^{\omega(c)}
\left(
\frac{1}{c} + \frac{1}{ s_q}
\right)
\end{equation}
Now,
$$
\sum_{\substack{c|d \\ c \ll s^{\cboundexponent}}}
\frac{s_c^{k-1}c^{1/2-o(1)}
D^{\omega(c)}}
{c}
\ll
\sum_{\substack{c|d \\ c \ll s^{\cboundexponent}}}
c^{-1/2+o(1)}
\ll
s_q^{o(1)}
$$
and similarly
$$
\frac{1}{ s_q}
\sum_{\substack{c|d \\ c \ll s^{\cboundexponent}}}
s_c^{k-1}c^{1/2-o(1)}
D^{\omega(c)}
\ll
\frac{1}{ s_q}
\sum_{\substack{c|d \\ c \ll s^{\cboundexponent}}}
c^{1/2+o(1)}
$$
Thus~(\ref{poly-eq:2}) is
\begin{multline}
\label{poly-eq:8}
\ll
\frac{s_q^{k-1}}{s_q^{k-1}}
\sum_{\substack{ d|q \\ d>s_q^T}}
\left(
\frac{s_q^{o(1)}}{d^{1/2-o(1)}}
+
\frac{1}{ s_q d^{1/2-o(1)}}
\sum_{\substack{c|d \\ c \ll s^{\cboundexponent}}}
c^{1/2+o(1)}
\right)
\\ =
s_q^{o(1)}
\sum_{\substack{ d|q \\ d>s_q^T}}
d^{-1/2+o(1)}
+
\frac{1}{s_q}\sum_{\substack{ d|q \\ d>s_q^T}}
\frac{1}{ d^{1/2-o(1)}}
\sum_{\substack{c|d \\ c \ll s^{\cboundexponent}}}
c^{1/2+o(1)}
\end{multline}
Now, for any $\beta \in (0,1/2)$,
\begin{multline*}
\sum_{\substack{ d|q \\ d>s_q^T}}
d^{-1/2+o(1)}
\ll
\sum_{d|q }
d^{-1/2+o(1)} \left( \frac{d}{s_q^T}  \right)^\beta
\\
\ll
s_q^{-\beta T}
\sum_{d|q }
d^{\beta-1/2+o(1)} 
\ll
s_q^{-\beta T+o(1)}.
\end{multline*}
Similarly, for any $\gamma > 0$, 
$$
\sum_{\substack{c|d \\ c \ll s^{\cboundexponent}}}
c^{1/2+o(1)}
\ll
s_q^{\gamma \cboundexponent}
\sum_{c|d}
c^{1/2-\gamma+o(1)}
\ll
s_q^{\gamma \cboundexponent}
d^{1/2-\gamma+o(1)}
$$
and thus 
$$
\sum_{\substack{ d|q \\ d>s_q^T}}
\frac{1}{ d^{1/2-o(1)}}
\sum_{\substack{c|d \\ c \ll s^{\cboundexponent}}}
c^{1/2+o(1)}
\ll
s_q^{\gamma \cboundexponent}
\sum_{ d|q }
d^{-\gamma+o(1)}
\ll
s_q^{\gamma \cboundexponent+o(1)}
$$
Hence (\ref{poly-eq:8}) is 
$$
\ll
s_q^{-\beta T+o(1)}
+
s_q^{-1 + \gamma \cboundexponent+o(1)}
\ll s_q^{-1/2+o(1)}
$$
if we take $T=1-o(1)$, 
$\beta = 1/(2T)-o(1)$, and $\gamma = 1/(2 \cboundexponent)$.
Thus, with $\alpha = 1/2 + o(1)$ (to bound the contribution from small
$d$), we find that
$$
\text{Error}_2 = 
\frac{1}{s_q^{k-1}}
\sum_{\substack{ d|q \\ d>1}}
\sum_{\vh \in s_q X' \cap \Z^{k-1}}
\eps_k(\vh,d)
\ll s_q^{-1/2+o(1)}
$$
\end{proof}




\end{document}